\documentclass{svmult}
\usepackage{amssymb,amsmath}
\usepackage{mathptmx}
\usepackage{math}
\usepackage{graphicx,pst-plot,epsfig,pstricks}
\usepackage{algorithm2e}
\newcommand{\LastUpdate}{4 September 2007}
\newcommand{\sfb}[1]{\textsf{#1}}

\begin{document}
\title*{Semantic
distillation: a method  for clustering objects by their contextual specificity}
\titlerunning{Last updated on \LastUpdate}
\author{Thomas Sierocinski\inst{1} \and Antony Le B{\'e}chec\inst{2} \and 
Nathalie Th\'eret\inst{2}  \and Dimitri Petritis\inst{1}
}
\authorrunning{Th.\ Sierocinski et al.}
\titlerunning{Semantic distillation}
\institute{Institut de recherche math\'ematique (UMR6625), Universit\'e de Rennes 1
\and
INSERM U620, Universit\'e de Rennes 1
}
\maketitle

\begin{abstract}
Techniques for data-mining, latent semantic analysis, contextual search of databases, etc.\ have long ago been developed by computer scientists working on 
information retrieval (IR). Experimental scientists, from all
disciplines, having to analyse large collections of raw experimental data
(astronomical, physical, biological, etc.) have developed
powerful methods for their statistical analysis and for  
clustering, categorising,
and classifying objects. 
Finally, physicists have developed a theory of quantum measurement,
unifying the logical, algebraic, and probabilistic aspects of queries
into a single formalism.

The purpose of this paper is twofold:
first to show that when formulated at an abstract
level, problems from IR,  from statistical data analysis, and from
physical measurement theories are 
very similar and hence can profitably be cross-fertilised, and,
secondly, to propose a novel method of fuzzy hierarchical clustering,
termed \textit{semantic distillation} --- strongly inspired from the
theory of quantum measurement ---, we developed to analyse raw
data coming from various types of experiments on DNA arrays. We 
illustrate the method by analysing DNA arrays experiments and
clustering the  genes of the array according to their specificity.
\vskip3mm

\end{abstract}

\section{Introduction}

Sequencing the genome constituted a culminating point in the analytic
approach of Biology. \cut{Its epistemological pendant in Physics is that
the fundamental constituents of matter are elementary particles.}
Now starts the era of the synthetic approach in Systems Biology where
interactions among genes induce their differential expression that leads
to the functional specificity of cells, the coherent organisation
of cells into tissues, organs, and finally organisms. 
\cut{The epistemological
pendant in Physics is the study of interactions among elementary particles
leading to the observed properties of matter.}

However, we are yet far from a complete
explanatory theory of living matter. It is therefore important to
establish precise and quantitative phenomenology before being able to
formulate a theory. The contribution of this paper is to provide the
reader with a novel algorithmic method, termed \textit{semantic distillation}, to analyse DNA arrays experiments 
(where genes are hybridised with various cell lines corresponding to various
tissues or specific individuals) by determining the degree of specificity of every gene  to the particular
context. The method 
provides experimental biologists with lists of candidate
genes (ordered by their degree of specificity) for every biological context,
clinicians with improved tools for diagnosis, 
pharmacologists with patient-tailored therapies, etc.

\cut{We propose in this paper a novel method of analysis, we term \textit{semantic
distillation}. Our method borrows ideas from several pre-existing algorithms;
it incorporates nevertheless novel features. 
It is furthermore very versatile, easily adaptable to various situations and
requirements.} In the sequel we present the method split into several 
algorithmic tasks thought as subroutines of the general algorithm. 
It is worth noting that the method, although can
profitably exploit, does not rely on any previous information stored
in the existing databases; its rationale is to help
analysing raw experimental data even in the absence of any previous
knowledge.

The main idea of the method is summarised as follows.
Experimental information hold on  the objects of the system undergoes a
sequence 
of processing steps; each step is performed on a different
representation of the information. 
Those different representation spaces and the corresponding
information processing act as successive filters revealing at the end
the most pertinent and significant part of the information, hence the name
``semantic distillation''.

At the first stage, raw experimental data, containing all available
information, are represented in an abstract Hilbert space, \textit {the
space of concepts} --- reminiscent of the space of pure states in 
Quantum Mechanics ---, endowing the set of objects with a metric space
structure that is exploited to quantify the interactions among objects and
encode them into a weighed graph on the vertex set of objects and with
object interactions as edge weights.

Now objects (genes) are parts of an organised system (cell, tissue, organism).
Therefore their mutual interactions are not just independent random variables;
they are interconnected through precise, although certainly very complicated
and mostly unknown  relationships.  
\cut{To give an image about what is meant here,
suppose that an observer knowing nothing about astronomy, measures
the positions in 3-dimensional space of a planet at $N$ different epochs, with 
$N$ being a large positive integer.
In the absence of any prior knowledge and measurement, 
the best prediction she can do 
is that the positions are $N$  \textit{random} 3-dimensional
vectors.
Now, positions of a planet are subject to the law 
of universal gravity attraction.   When our observer plots the $N$ experimental
3-dimensional measurements they appear organised on a 1-dimensional 
manifold, corresponding to the elliptic orbit of the planet around the sun.
The  (unknown to our observer) law of universal attraction
imposes relations among the observed positions that project the
distribution of points from an \textit{a priori} 3-dimensional random
cloud to an ellipse.}
 We seek to reveal (hidden and unknown) interactions
among genes. This is achieved by trading the weighed graph representation for
a low-dimensional representation and using spectral properties of the
weighed Laplacian on the graph to grasp the essential interactions.

The following step consists in a fuzzy divisive clustering of objects among two
subsets by exploiting the previous low-dimensional representation. This procedure
assigns a fuzzy membership to each object relative to characters
of the two subsets. Fuzziness is as a matter of fact a distinctive property
of experimental biological data reflecting our incomplete knowledge of
fundamental biological processes. 

Up to this step, our method is a sequence of known algorithms that have been
previously used separately in the literature in various contexts. The novelty
of our method relies on the following steps.
The previous fuzzy clustering reduced the indeterminacy of the system.
This information is fed back to the system to perform a projection to
a proper Hilbert subspace.
In that way, the information content of the dataset is modified by the information gained by the previous observations.  After this feeding back, the three
previous steps are repeated but now referring to a Hilbert spaces of
lower dimension. Therefore our method is not a mere fuzzy clustering
algorithm but a genuine non-classical interaction  information retrieval
procedure where previous observations alter the informational content
of the system, reminiscent of the measurement procedure in Quantum Mechanics.

\section{A Hilbert space formulation}

\subsection{Mathematical form of the dataset}

Let $\mathbb{B}$ be a finite set of \textit{documents}
(or objects, or books) and $\mathbb{A}$ a finite set of \textit{attributes}
(or contexts, or keywords).
The dataset is a $|\mathbb{B}|\times|\mathbb{A}|$ matrix
$\sfb{X}=(x_{ba})_{b\in\mathbb{B}, a\in\mathbb{A}}$
of real or complex elements, where $|\cdot|$ represents cardinality.
Equivalent ways of representing the dataset are 
\begin{itemize}
\item
a collection
of $|\mathbb{B}|$ row vectors 
$\bx_b=(x_{b1},\ldots,x_{b|\mathbb{A}|}), b\in\mathbb{B}$ of $\mathbb{R}^{|\mathbb{A}|}$ 
(or $\mathbb{C}^{|\mathbb{A}|}$),
\item
a collection
of $|\mathbb{A}|$ column vectors 
$\bx^a=(x_{1a},\ldots,x_{|\mathbb{B}|a})
, a\in\mathbb{A}$ of $\mathbb{R}^{|\mathbb{B}|}$ 
(or $\mathbb{C}^{|\mathbb{B}|}$).
\end{itemize}

\begin{example}
In the experiments we analysed
$\mathbb{B}$ is  a set of 12000 human genes and 
$\mathbb{A}$ a set of 12
tissular contexts. The matrix elements $x_{ba}$ are real
numbers encoding luminescence intensities (or their logarithms)
of DNA array ultimately
representing the level of expression of gene $b$ in context $a$.
\end{example}

\begin{example}
Let 
$\mathbb{B}$ be a set of books in a library and $\mathbb{A}$ a set of 
bibliographic keywords. 
The matrix elements $x_{ba}$ can be $\{0, 1\}$-valued:  if the term
$a$ is present in the book $b$ then 
$x_{ba}=1$ else  $x_{ba}=0$.
A variant of this example is when $x_{ba}$ are integer valued:
if the term $a$ appears $k$ times in document $b$ then $x_{ba}=k$.
\end{example}

\begin{example}\label{exam:students}
Let 
$\mathbb{B}$ be a set of students and $\mathbb{A}$ a set of 
papers they gave. 
The matrix elements $x_{ba}$ are real valued; $x_{ba}$ is the
mark the student $b$ got in paper $a$.
\end{example}

The previous examples demonstrate the versatility of the
method by keeping the formalism 
at an abstract level to apply
indistinctively into various very different situations without any change.
Note also that the assignment as set of documents or attributes is a matter
of point of view; for instance, example \ref{exam:students} as it stands
is convenient in evaluating students. Interchanging the role of
sets $\mathbb{A}$ and $\mathbb{B}$ renders it adapted to the evaluation
of teaching. As a rule of thumb, in biological applications,  
 $|\mathbb{A}|\ll|\mathbb{B}|$.

\subsection{The space of concepts}

For $\mathbb{A}$
and $\mathbb{B}$ as in the previous subsection, we define 
the \textit{space of concepts}, $\mathcal{H}_\mathbb{A}$, as the real
or complex free vector space over $\mathbb{A}$, i.e.\ elements
of $\mathbb{A}$ serve as indices 
of an orthonormal basis of $\mathcal{H}_\mathbb{A}$.
Therefore, the complete dataset $\sfb{X}$ can be represented
as the collection of $|\mathbb{B}|$ vectors
$\lar{\Xi_b}=\sum_{a\in\mathbb{A}} x_{ba} \lar{a}\in 
\mathcal{H}_\mathbb{A}$, with
$b\in \mathbb{B}$ and where $\lar{a}$ represents
the element of the orthonormal basis of the free vector space
corresponding to the attribute $a$.
We use here  Dirac's notation to represent vectors, 
linear forms and projectors
on this space (see any book on quantum mechanics or \cite{petritis-iq} 
for a freely accessible document and  \cite{vanRijsbergen} for
the use of this notation in information retrieval).
The vector $\lar{\Xi_b}$ contains all available experimental information
on document $b$ in various cellular contexts indexed by the attributes
$a$; it can be thought as a convenient bookkeeping device of the data
$(x_{ba})_{a\in \mathbb{A}}$, in the same way a generating function contains
all the information on a sequence as formal power series.
 
The vector space is equipped with a scalar product
defined for every two vectors 
$\lar{\psi}= \sum_{a\in\mathbb{A}} \psi_{a} \lar{a}$
and 
$\lar{\psi'}= \sum_{a\in\mathbb{A}} \psi'_{a} \lar{a}$
by 
$\scalar{\psi}{\psi'}=\sum_{a\in\mathbb{A}} \overline{\psi}_a \psi'_{a}$,
where $ \overline{\psi}_a$ denotes the complex conjugate of $\psi_a$
(it coincides with $\psi_a$ if it is real).
Equipped with this scalar product, the vector space 
$\mathcal{H}_\mathbb{A}$ becomes a real or complex
$|\mathbb{A}|$-dimensional Hilbert space.
The scalar product induces a Hilbert norm on the space, denoted by $\|\cdot\|$. 
In the sequel
we introduce also \textit{rays} on the Hilbert space i.e.\ 
normalised vectors. Since the dataset $\sfb{X}$ does not in principle
verify any particular numerical constraints,
rays are constructed by dividing vectors by their norms.
We use the symbol $\lar{\xi_b}=\lar{\Xi_b}/\|\lar{\Xi_b}\|$ to denote
the ray associated with vector $\lar{\Xi_b}$.

The Hilbert space structure on   $\mathcal{H}_\mathbb{A}$
allows a natural geometrisation of the space of
documents by equipping it with a pseudo-distance\footnote{It 
is termed pseudo-distance since it verifies
symmetry and triangle inequality but  $d(b,b')$
can vanish
even for different $b$ and $b'$. 
As a matter of fact, $d$ is a distance on
the projective Hilbert space.}
$d:\mathbb{B}\times\mathbb{B}\rightarrow \mathbb{R}_+$
defined by $d(b,b')=\|\lar{\xi_b}-\lar{\xi_{b'}}\|$.
What is important here is not the precise form of the pseudo-metric
structure of $(\mathbb{B},d)$; several other pseudo-distances
can be introduced, not necessarily compatible with the scalar product.
In this paper we stick however to the previous pseudo-distance,
postponing into a later publication explanations about the 
significance of other pseudo-distances.

As is the case in Quantum Mechanics, 
the Hilbert space description incorporates into a unified algebraic
framework all logical and probabilistic information hold by the
dataset. An enquiry of the type ``does the system possess feature $F$''
is encoded into a projector $P_F$ acting on the Hilbert space. The subspace
associated with the projector $P_F$ is interpreted as the set of
documents retrieved by asking the question about the feature $F$. Now
all experimental information hold by the dataset is encoded into the
\textit{state} of the system represented by a 
\textit{density matrix} 
$\rho$ (i.e.\ a self-adjoint, positive, trace class
operator acting on $\mathcal{H}_{\mathbb{A}}$ having  unit trace).
Retrieved documents possess the feature $F$ with probability
$\tr (\rho P_F)$. Thus the algebraic description incorporates 
logical information on the documents retrieved as relevant to a given
feature and assign them a probability determined by the state defined
by the experiment. For example, the probability that a gene $b$
is relevant to an attribute $a$ is given by the above formula with
$P=\lar{a}\sca{a}$ and $\rho=\lar{\xi_b}\sca{\xi_b}$, yielding
$\tr(\rho P)=|\scalar{\xi_b}{a}|^2$.

\section{A weighed graph with augmented vertex set}
The careful reader has certainly already noted
that in the above description vectors $\lar{\xi_b}$, encoding
the information about document $b$, and 
basis vectors $\lar{a}$, associated with attribute $a$,
all belong to the same Hilbert space $\mathcal{H}_\mathbb{A}$.
Therefore, although initially the sets $\mathbb{A}$ and $\mathbb{B}$
are disjoint since they 
have distinct elements, 
when passing to the Hilbert space representation, vectors
$\lar{\xi_b}$ and $\lar{a}$ have very similar roles in representing
indistinguishably objects or 
attributes as vectors of $\mathcal{H}_\mathbb {A}$.
In the sequel, we introduce the set $\mathbb{V}$ (or more precisely
$\mathbb{V}_\mathbb{A}$ to remove any ambiguity) as the set
$\mathbb{V}_\mathbb{A}=\mathbb{A}\cup\mathbb{B}.$
Thus, for any $v\in\mathbb{V}_\mathbb{A}$,
\[\lar{\Xi_v}=\left\{\begin{array}{ll}
\lar{a} & \textrm{if }\ 
v=a\in \mathbb{A},\\
\sum_{a\in\mathbb{A}} x_{ba}\lar{a}
& \textrm{if }\ 
v=b\in \mathbb{B}.
\end{array}
\right.\]
The new vectors $\lar{\Xi_a}=\lar{a}$
are included as \textit{specificity witnesses}
in the dataset. Note that since these new vectors
are also elements of the same Hilbert space, the pseudo-distance
$d$ naturally extends to $\mathbb{V}_\mathbb{A}$.

Suppose now that a \textit{similarity function} $\sigma:
\mathbb{V}_\mathbb{A}\times
\mathbb{V}_\mathbb{A}\rightarrow [0,1]$ is defined.
For the sake of definiteness, the reader can think of $\sigma$
as being given, for example\footnote{This function is well adapted to datasets
$\sfb{X}=(x_{ba})$, with $x_{ba}\in \mathbb{R}^+$; for more general
datasets, the factor $1/2$ must be changed to $1/4$.}, 
by $\sigma(v,v')=\sqrt{1-\frac{1}{2}
d(v,v')^2}$; results we quote in section \ref{sec:results}
are obtained with a slight modification of this similarity function.  
However, again, the precise
form of the similarity function is irrelevant in the abstract setting
serving as foundation of the method. Several other similarity functions
have been used like, for example, 
$\sigma(v,v')=\exp(-\|\Xi_v-\Xi_{v'}\|^2/\tau)$ with $\tau$ a
positive constant or some others,  
in particular, functions taking value 0 even for
some vertices corresponding to non orthogonal rays but the explanation
of their
significance is postponed to a subsequent publication.

A weighed graph is now constructed with vertex set 
$\mathbb{V}_\mathbb{A}$. Weights are assigned to the edges
of the complete graph over $\mathbb{V}_\mathbb{A}$; the weights
being expressible in terms of the similarity function $\sigma$.
Again, the precise expression is irrelevant for the exposition of the
method. For the sake of concreteness, the reader can suppose
that the weights $W_{vv'}$ are given by
$W_{vv'}=\sigma(v,v')$. 
The pair $(\mathbb{V}_\mathbb{A},W)$ with $W$ being the
symmetric matrix $W=(W_{vv'})_{v,v\in \mathbb{V}_\mathbb{A}}$,
denotes the weighed graph.

At this level of the description we follow now standard techniques
of reduction of the data dimensionality by optimal representation
of the graph in low dimensional Euclidean spaces spanned
by eigenvectors of the Laplacian. Such methods
have been used by several authors \cite{Belkin,Nilsson}.
Here we give only the basic definitions and main results of this method. 
The interested reader may consult standard textbooks
like \cite{Chung,Cvetkovic,Godsil} for general exposition of the
method.

\begin{definition}
A map $\br:\mathbb{V}_\mathbb{A}\rightarrow \mathbb{R}^\nu$
is called a \textit{$\nu$-dimensional representation} of the graph. 
The representation is always supposed non-trivial
(i.e.\ $\br\not\equiv 0$)
and balanced (i.e.\ $\sum_{v\in\mathbb{V}_\mathbb{A}} \br(v)=0$).
\end{definition}
From the  weights matrix $W$ we construct the \textit{weighed Laplacian matrix} $\Lambda=D-W$ where 
the matrix elements $D_{vv'}$ are 0 if $v\ne v'$ and
equal to $\sum_{v''\in \mathbb{V}_\mathbb{A}} W_{vv''}$ if
$v=v'$.
More precisely, we denote by $\Lambda(\mathbb{V}_\mathbb{A})$ this
weighed Laplacian to indicate that it is defined on the vertex set 
$\mathbb{V}_\mathbb{A}$. This precision will be necessary in the next
section specifying the semantic distillation algorithm where
the vertex set will be recursively modified at each step.
The \textit{weighed energy of the representation} is given by
\[\mathcal{E}_W(\br)=\sum_{v,v'\in \mathbb{V}_\mathbb{A}}
W_{vv'} \|\br(v)-\br(v')\|^2,\]
where in this formula
$\|\cdot\|$ denotes the Euclidean norm of $\mathbb{R}^\nu$.

\begin{theorem}
Let $N=|\mathbb{V}_\mathbb{A}|$ and $\{\lambda_1,
\ldots, \lambda_N\}$ be the spectrum of $\Lambda$, ordered
as $\lambda_1\leq\lambda_2\ldots\lambda_N$. Suppose that
$\lambda_2>0$. Then 
$\inf_\br\mathcal{E}_W(\br)=\sum_{i=2}^{\nu+1}\lambda_i$,
where the infimum is over all $\nu$-dimensional non-trivial balanced
representations of the graph.
\end{theorem}
\begin{remark} If $\bu^1, \ldots
\bu^N$ are the eigenvectors of $\Lambda$ corresponding to the
eigenvalues $\lambda_1,\ldots,\lambda_N$ ordered as above,
then $\bu^2$ is the best one-dimensional,
$[\bu^2,\bu^3]$ the best two-dimensional, etc,
$[\bu^2,\ldots, \bu^{\nu+1}]$
the best $\nu$-dimensional representation of the graph
$(\mathbb{V}_\mathbb{A}, W)$.
\end{remark}

\section{Fuzzy semantic clustering and distillation}
\label{sec:semantic-distillation}

The algorithm of semantic distillation is 
a recursive divisive fuzzy clustering followed
by a projection on a Hilbert
subspace and a thinning of the graph. 
It starts with the Hilbert space $\mathcal{H}_{\mathbb{A}}$
and the graph with vertex set $\mathbb{V}_{\mathbb{A}}$
and constructs
a sequence of Hilbert subspaces and subgraphs
indexed by the words $\kappa$ of finite length
on a two-letter alphabet. This set is isomorphic to a subset of the rooted binary tree.
If $\kappa$ is the root, then define $\mathbb{M}_\kappa=\mathbb{A}$.
Otherwise, $\mathbb{M}_\kappa$ will be a proper
subset of $\mathbb{A}$, i.e.\ $\emptyset\subset\mathbb{M}_\kappa
\subset \mathbb{A}$,indexed by $\kappa$.
When $|\mathbb{M}_\kappa|=1$ then the corresponding $\kappa$
is a leaf of the binary tree. The algorithm stops when all indices
correspond to leaves. 

More precisely, let $\mathbb{K}=\{1,2\}$, 
$\mathbb{K}^0=\{\kappa:\kappa=()\}$,
and for integers $n\geq 1$ let 
$\mathbb{K}^n=\{\kappa:\kappa=\kappa_1\cdots\kappa_n; 
\kappa_i\in\mathbb{K}\}$. Finally let 
$\mathbb{K}^*=\cup_{n\geq 0} \mathbb{K}^n$ denote the set of
words on two letters of indefinite length, 
including the empty sequence, denoted by $()$, of zero length that coincides with the root
of the tree. If $\kappa=\kappa_1\cdots\kappa_n$ is a word of $n$ letters
and $k\in\mathbb{K}$, we denote the concatenation $\kappa k$ as the word
of $n+1$ letters $\kappa_1\cdots\kappa_nk$.

We start from the empty set $\textsf{Leaves}=\{\}$, the empty sequence
$\kappa=()$ and the current attributes set $\mathbb{M}_\kappa=
\mathbb{M}_{()}=\mathbb{A}$ and current tree $\textsf{Tree}=\{\kappa\}$. 
We denote $\mathbb{V}_\kappa=
\mathbb{B}\cup\mathbb{M}_\kappa$. We need further
a \textit{fuzzy membership} function $m:\mathbb{V}_\kappa\times\mathbb{K}
\rightarrow [0,1]$. The fuzzy clustering algorithm is succinctly described
as Algorithm \ref{alg:FC} below.

\begin{algorithm}\caption{\label{alg:FC}\textsf{FuzzyClustering}}
\SetLine
\KwData{$\kappa$, $\mathbb{M}_\kappa$, $\br$, 
objective function $F$ }
\KwResult {Two sets $\mathbb{M}_{\kappa1}$ and $\mathbb{M}_{\kappa2}$
and the fuzzy membership  $m(v,k)$ for $v\in \mathbb{V}_\kappa$
in the clusters $\mathbb{M}_{\kappa1}$ and $\mathbb{M}_{\kappa2}$}
\If {$|\mathbb{M}_\kappa|>1$}
{
{\textbf{assign}
$(v_1,v_2)\leftarrow \arg\max\{\|\br(v)-\br(v')\|,
v,v'\in \mathbb{V}_\kappa\}$}\;
{\textbf{assign} $\br(v_1)$ and $\br(v_2)$ as centroids
for the two candidate finer clusters 
$\mathbb{M}_{\kappa1}$ and $\mathbb{M}_{\kappa2}$}\;
{\textbf{use} standard $2$-means fuzzy clustering algorithm  
to minimise objective function $F$ under the constraint $\sum_{k=1}^2
m(v,k)=1$, for all $v\in\mathbb{M}_\kappa$ }\;
{\textbf{assign} $\mathbb{M}_{\kappa1}\leftarrow\{v\in \mathbb{M}_\kappa:
m(v,1)>m(v,2)\}$}\;
{\textbf{assign} $\mathbb{M}_{\kappa2}\leftarrow
\mathbb{M}_\kappa\setminus \mathbb{M}_{\kappa1}$}\;
 }
\end{algorithm}
\vskip5mm

Note that in the previous construction $\mathbb{M}_{\kappa k}\subset
\mathbb{M}_{\kappa}$ for every $\kappa$ and every $k\in\mathbb{K}$.
Therefore, the algorithm explores the branches of a tree from the root
to the leaves.
Denote by $\pi_\kappa$ the orthogonal projection from
$\mathcal{H}_\mathbb{A}$ to $\mathcal{H}_{\mathbb{M}_\kappa}$.
The distillation step is described by the following Algorithm \ref{alg:SD}.

\begin{algorithm}\caption{\label{alg:SD}\textsf{Distillation}}
\SetLine
\KwData {\textsf{FuzzyClustering} }
\KwResult {\textsf{Leaves} and sequence of singleton sets
$\mathbb{M}_{\kappa}$ for $\kappa\in\textsf{Leaves}$}

\textbf{Initialisation}\{\\
$\kappa\leftarrow ()$\;
$\mathbb{M}_\kappa\leftarrow \mathbb{A}$\;
$\textsf{Leaf}(\kappa)\leftarrow\mathbb{M}\kappa$\;
$\textsf{Leaves}\leftarrow \{\}$\;
$\textsf{Tree}\leftarrow \{\kappa\}$\;
$\textsf{Bookkeeping}\leftarrow \{\kappa\}$\;
\}\\
\While{$\textsf{Bookkeeping}\ne\emptyset$}{
    \For{$\kappa\in\textsf{Bookkeeping}$}{
        \eIf{ $|\mathbb{M}_\kappa|=1$}{
             $\textsf{Leaves}\leftarrow\textsf{Leaves}\cup\{\kappa\}$\;
             $\textsf{Bookkeeping}\leftarrow\textsf{Bookkeeping}\setminus\{\kappa\}$\;}           {
\textbf {Use} $\pi_\kappa$ to project from 
$\mathcal{H}_\mathbb{A}$ to $\mathcal{H}_{\mathbb{M}_\kappa}$\;
\textbf{Thin} the graph:
 $\mathbb{V}_\kappa\leftarrow\mathbb{B}\cup\mathbb{M}_\kappa$\; 
\textbf {Compute} weighed Laplacian $\Lambda(\mathbb{V}_\kappa)$\;
\textbf{Diagonalise} $\Lambda(\mathbb{V}_\kappa)$\;
 \textbf{Compute}
$\nu$-dimensional representation $\br$\;
 \textbf{Call}            \textsf{FuzzyClustering}\;
             \For{$k\in\mathbb{K}$}{
                  $\kappa'\leftarrow\kappa k$\;
                  $\textsf{Leaf}(\kappa')\leftarrow \mathbb{M}_{\kappa'}$ /* 
$\mathbb{M}_\kappa$ as determined by \textsf{FuzzyClustering} */\;
$\textsf{Tree}\leftarrow\textsf{Tree}\cup\{\kappa'\}$\;
$\textsf{Bookkeeping}\leftarrow\textsf{Bookkeeping}\cup\{\kappa'\}$\;
}}}}
\end{algorithm}
\vskip5mm

\cut{
To illustrate the method,
suppose thus that we are given the eigenvector $u^2$
corresponding to the second smallest eigenvalue of $\Lambda$,
providing us with the best one-dimensional representation of the
graph. Now, $u^2$ is a vector of $\mathbb{R}^N$, with 
$N=|\mathbb{D}_\mathbb{A}|$.
Suppose we order the $N$ components of $u^2$ in ascending
order. This ordering corresponds to a bijective relabelling 
$r:\mathbb{D}_\mathbb{A}\rightarrow \{1, \ldots, N \}$
of the vertices so that for any $i,j\in \{1, \ldots, N \}$
with $i<j$, we have $u^2(r^{-1}(i))<u^2(r^{-1}(j))$.
Therefore the set of vertices becomes now a totally ordered set
isomorphic to the set $\{1, \ldots, N \}$.

Now recall that the set $\mathbb{D}_\mathbb{A}$ is the set
$\mathbb{D}$ of documents, augmented by the set
$\mathbb{A}$ of attributes, the latter acting as 
specificity witnesses. In the relabelling induced by $u^2$, 
the specificity witnesses are relabelled as elements
$\{r(a_1),\ldots, r(a_{|\mathbb{A}|})\}=
\{\alpha_1,\ldots,\alpha_{|\mathbb{A}|}\}
\subseteq \{1, \ldots, N \}$, with $\alpha_1<\ldots<\alpha_{|\mathbb{A}|}$.
Let 
\[i=\arg\max\{u^2(r^{-1}(\alpha_j))-u^2(r^{-1}(\alpha_{j-1})): 
j=2,\ldots, |\mathbb {A}|\}\]
correspond to the most frank separation of the attributes and define
$\mathbb{A}_1=\{r^{-1}(\alpha_1),\ldots,
r^{-1}(\alpha_{i-1})\}$ and $\mathbb{A}_2=\mathbb{A}\setminus
\mathbb{A}_2$.  
Cluster the documents by assigning
\[C_{\mathbb{A}_1}=
\{d\in\mathbb{D}: \forall a \in \mathbb{A}_2,
|u^2(d)-u^2(a)|>\min_{a'\in\mathbb{A}_1}|u^2(d)-u^2(a')|\},\]
and similarly for $C_{\mathbb{A}_2}$ by exchanging the roles of $\mathbb{A}_1$ 
and $\mathbb{A}_2$.
The set $C_{\mathbb{A}_1}$ contains the documents 
the most specific to attributes $\mathbb{A}_1$ and similarly for the
set $C_{\mathbb{A}_2}$.  

Semantic distillation is a recursive procedure 
producing partitions $\mathbb{A}_1$ and $\mathbb{A}_2$
of subsets of $\mathbb{A}$ repeated until
both sets $\mathbb{A}_1$ and $\mathbb{A}_2$ become singletons.
The first step of semantic distillation corresponds 
to an orthogonal projection $P=\sum_{a\in\mathbb{B}}\larsca{a}{a}$ 
from the 
Hilbert space $\mathcal{H}_{\mathbb{A}}$ onto
the Hilbert subspace $\mathcal{H}_{\mathbb{B}}$, where
the set $\mathbb{B}$ is the smallest of the above defined sets
$\mathbb{A}_1$ and $\mathbb{A}_2$. 
The graph construction is started afresh 
on the vertex set of the \textit{thinned graph}
$(\mathbb{D}\setminus C_{\mathbb{B}})_{\mathbb{B}^c}$ 
and weights computed in terms of the Hilbert norm on the Hilbert subspace
$\mathcal{H}_\mathbb{B}$.
This gives rise to a new weight matrix $W$ and a new Laplacian 
$\Lambda\left((\mathbb{D}\setminus C_{\mathbb{B}})_{\mathbb{B}^c}\right)$.
We proceed then with the computation of eigenvalues and eigenvectors
of this new Laplacian leading to separation of remaining attributes
into two new subsets and new clustering. The method ends
when all objects are clustered.
}
 
\section{Illustration of the method, robustness and complexity issues}
\label{sec:results}
We tested the method on a dataset for an experiment on DNA array
published in \cite{Yanai},  with   the set
$\mathbb{A}$ of attributes corresponding to 12 cell lines 
(bone marrow, liver, heart, spleen, 
lung, kidney, skeletal muscle, 
spinal cord, thymus, brain, prostate, pancreas) and
the set $\mathbb{B}$ of documents corresponding
to 12000 human genes. 
\cut{In the figure \ref{fig:relabelling}, we present two different steps of
the eigenvector components corresponding to best one-dimensional
representations (i.e.\ $\nu=1$) for the thinned graph. 
Since the representation is $1$-dimensional, this ordering
corresponds to a relabelling of genes.

\begin{figure}
\begin{center}
 \begin{tabular}{c c}
   \rotatebox{270}{\includegraphics[height=5cm,width=3cm]{fig-IR-liver_temp.eps}} &
   \rotatebox{270}{\includegraphics[height=5cm,width=3cm]{fig-IR-skeletalmuscle_temp.eps}}\\
                \small{liver} & \small{skeletal muscle}\\
        \end{tabular}
\end{center}
\caption{\label{fig:relabelling}
Abscissas represent a relabelling of genes, ordinates
the component of $\bu^2$ vector on the corresponding genes. The cells 
indicated under the figures
represent the attributes associated with the smallest of the sets
$\mathbb{A}_1$ or $\mathbb{A}_2$ at the considered intermediate step. It turns
out that these smallest intermediate sets are singletons in both shown cases.}
\end{figure}
}
To illustrate the method we present here only an example of the type of results
we obtain by our method for the simplest case of one-dimensional representation
of the graph. The complete lists of specificity degrees for the
various genes (including their UniGene identifiers) for various dimensions
are provided as supplemental material (at the home page of the first author).

\begin{figure}
\begin{center}
 \begin{tabular}{c c}
   \rotatebox{270}{\includegraphics[height=5cm,width=3cm]{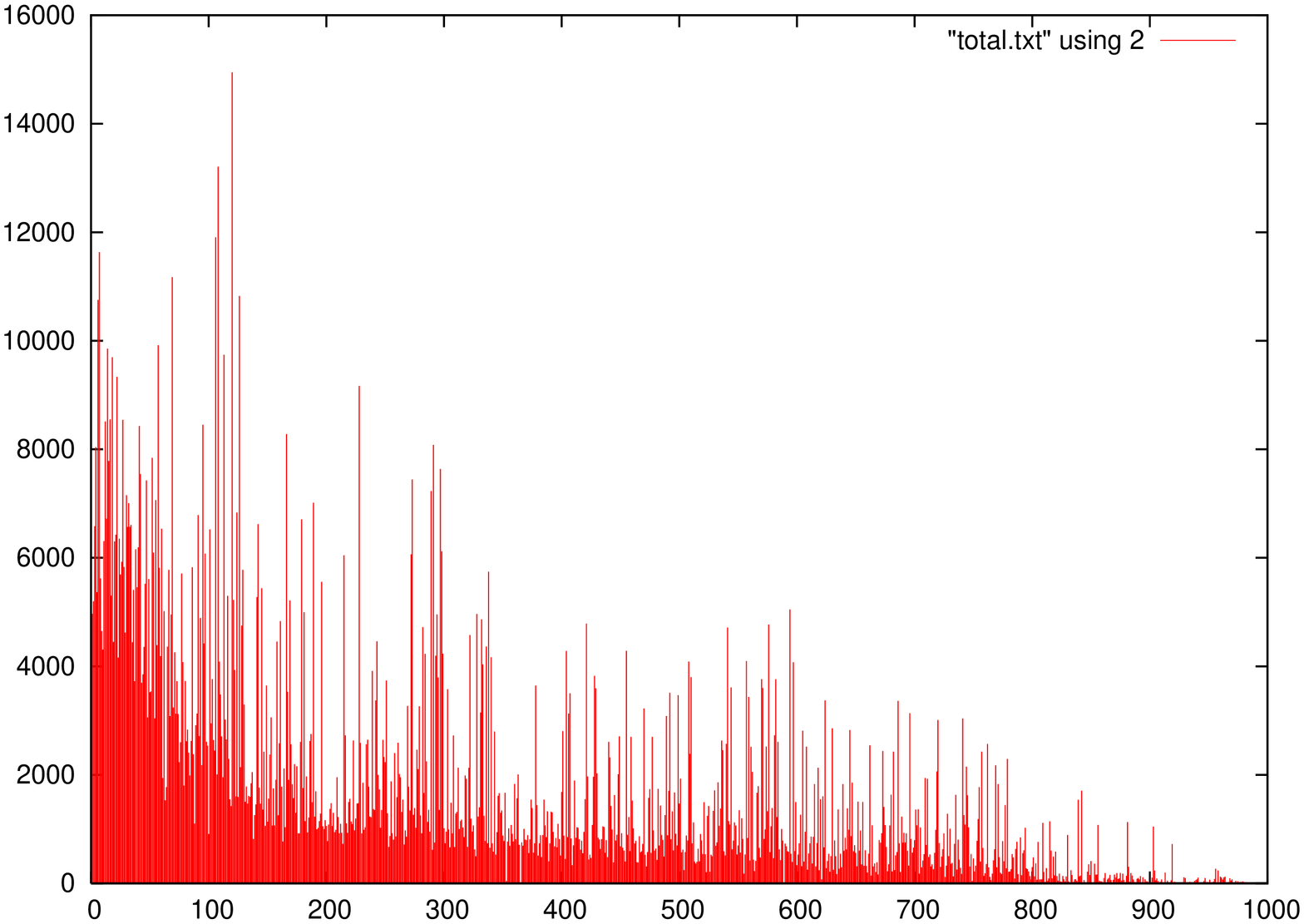}} &
   \rotatebox{270}{\includegraphics[height=5cm,width=3cm]{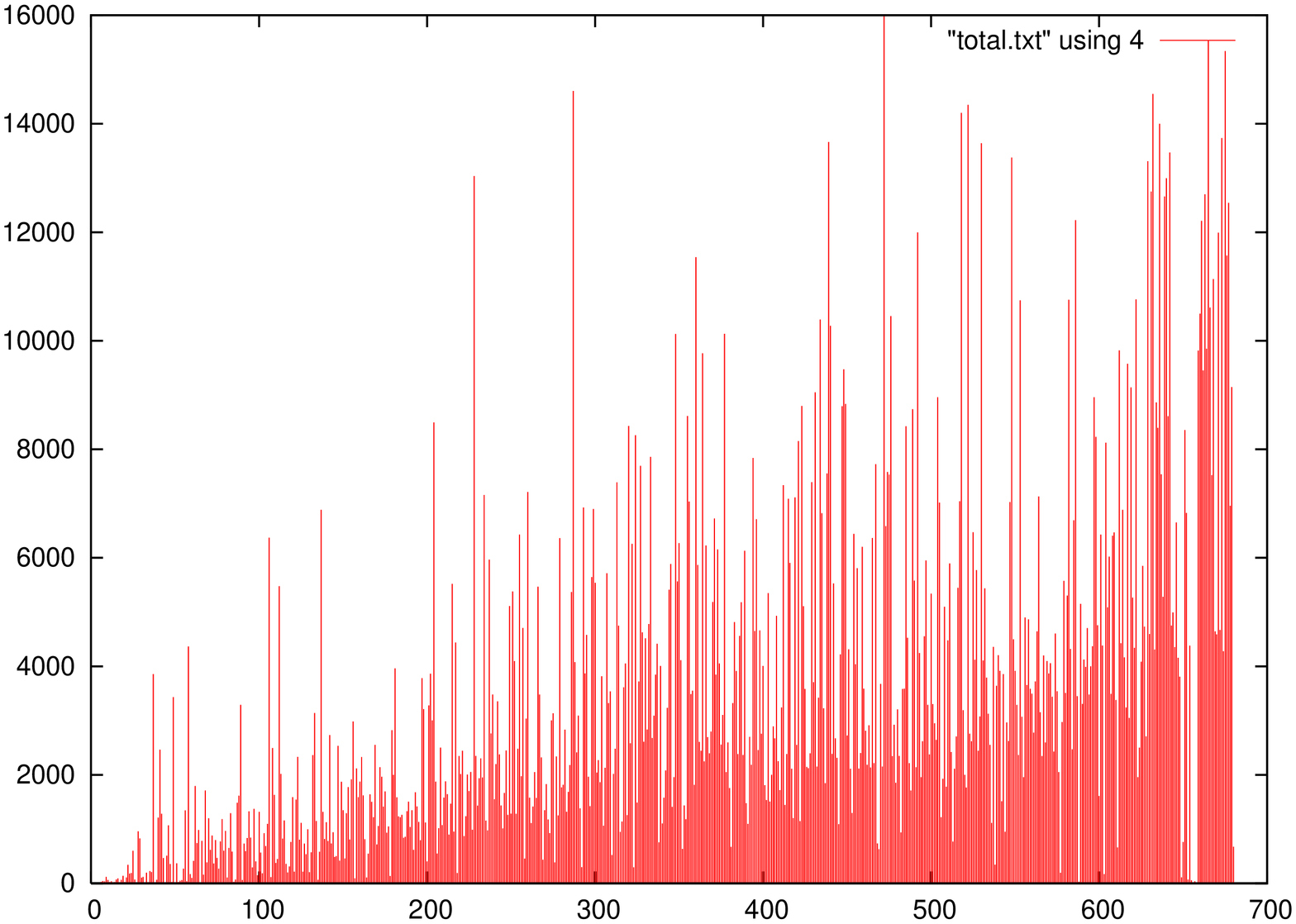}}\\
                \small{liver} & \small{skeletal muscle}\\
        \end{tabular}
\end{center}
\caption{\label{fig:clustering}
 For every singleton cluster, i.e.\ tissular context $\kappa\in\textsf{Leaves}$
(we  present solely the cases
$\mathbb{M}_\kappa=\{\textrm{liver}\}$ 
and $\mathbb{M}_\kappa=\{\textrm{ skeletal muscle}\}$ in this example), the 
horizontal axis contains the set $\mathbb{B}$ of genes \textit{relabelled}
according to their decreasing (resp. increasing) fuzzy membership to 
$\mathbb{M}_\kappa$.
Vertical axis represents the 
experimentally measured level of expression for those genes.
}
\end{figure}

Note that for one-dimensional representation, ordering by the magnitude of the
eigenvector components is equivalent to
a relabelling of genes.
The figure \ref{fig:clustering} represents, within the previous
mentioned relabelling, the
levels of expressions for clustered genes. 
The same procedure has been applied for higher dimensional representation
of the graph (i.e.\ $\nu>1$). These results are not presented here; they 
marginally improved some specifications and helped us removing apparent
degeneracy in some cases.
Finally, in the table \ref{tab:annotation},
we give an example of the annotation provided by the database UniGene
for the genes classified as specific of skeletal muscle cell line by our method.
\begin{table}
\caption{\label{tab:annotation} Annotation of the genes
closest (within the relabelling induced by $\bu^2$)
to the specificity witness ``skeletal muscle''. Genes are separated by the
-- symbol.}
\hrule
\vskip2mm
{\scriptsize
\textit{ATPase, Ca++ transporting, cardiac muscle, fast twitch 1, calcium signaling pathway -- Troponin I type 2; skeletal, fast -- Myosin, light chain 1, alkali; skeletal, fast -- Ryanodine receptor 1; skeletal; calcium signaling pathway -- Fructose-1,6-bisphosphatase 2, glycolysis / gluconeogenesis -- Actinin, alpha 3; focal adhesion -- Adenosine monophosphate deaminase 1 (isoform M) purine metabolism -- Troponin C type 2; fast; calcium signaling pathway -- Carbonic anhydrase III, muscle specific; nitrogen metabolism -- Nebulin -- Troponin I type 1; skeletal, slow -- Myosin, heavy chain 3, skeletal muscle -- Myogenic factor 6, herculin -- Myosin binding protein C, fast type -- Calcium channel, voltage-dependent, beta 1 subunit -- Metallothionein 1X -- Bridging integrator 1 -- Bridging integrator 1 -- Calpain 3, (p94) -- Tropomyosin 3 -- Phosphorylase, glycogen; muscle (McArdle syndrome, glycogen storage disease type V); starch and sucrose metabolism -- Myozenin 3 -- Myosin binding protein C, slow type -- Troponin T type 3; skeletal, fast -- Superoxide dismutase 2;  mitochondrial -- Nicotinamide N-methyltransferase -- Sarcolipin -- Interleukin 32 -- Sodium channel, voltage-gated, type IV, alpha subunit -- Guanidinoacetate N-methyltransferase; urea cycle and metabolism of amino groups.}}\\
\hrule
\end{table}

We observe that the majority of 
genes classified as most specific by our method
are in fact annotated as specific in the database. 
To underline the power of our method, note that the UniGene annotation for
the ATPase gene is ``cardiac muscle''. Our method determines it as 
most specific of ``skeletal
muscle''. We checked the experimental data we worked on and realised that 
this gene is, as a matter of fact,  
5 times more expressed in the skeletal muscle context than
in the cardiac muscle. Therefore, our method correctly determines this
gene as skeletal-muscle-specific. 

In summarising, our method
is an automatic and algorithmic method of analysis of raw experimental data; it can be used to any experiment
of similar type \textit{independently of any previous knowledge included
in genomic databases} to provide biologists
with a powerful tool of analysis. In particular, since most
of the genes are not yet annotated in the existing databases, the method 
provides biologists
with candidate genes for every particular context for further investigation. 
Moreover, the genetic character of documents and attributes
is purely irrelevant; the same method can be used to any other
dataset of similar structure, let them concern linguistic, genetic, or image
data. 

Concerning the algorithmic complexity 
of the method, the dominant contribution 
comes from the diagonalisation of a $|\mathbb{B}|\times
|\mathbb{B}|$ dense real symmetric matrix, requiring at worst
$\mathcal{O}(|\mathbb{B}|^3)$ time 
steps and $\mathcal{O}(|\mathbb{B}|^2)$ space.
The time complexity can be slightly reduced, if only low-dimensional 
(dimension $\nu$) 
representations are sought,  
to $\mathcal{O}(\nu\times|\mathbb{B}|^2)$ time steps.
Moreover, we tested the method against additive or multiplicative random 
perturbations of the experimental data; it proved astonishingly robust.

\section{Connections to previous work}

The algorithm of semantic distillation maps the
dataset into a graph and uses spectral methods and fuzzy clustering 
to analyse the
graph properties. As such, this algorithm is inspired by various pre-existing
algorithms and borrows several elements from them. 

The oldest implicit use of a vector space structure to represent dataset
and application of spectral methods to analyse them is certainly ``principal
components analysis'' introduced in \cite{Pearson}. The method seeks
finding 
directions of maximal variability in the space corresponding to
linear combinations of the underlying vectors.\cut{; an expository account
of the method can be found in \cite{Anderson}. }
The major drawbacks of principal components analysis are
the assumptions that dataset matrix is composed of row vectors
that  are independent and identically distributed realisations of
the same random vector (hence the covariance matrix whose principal components
are sought can be approximated by the empirical covariance of the process)
and that there exists a linear transformation maximising the variability.

Vector space representations and singular value
decomposition,  as reviewed in \cite{Berry},
have been used to retrieve information from digital libraries. 
Implementations
of these ideas range from the famous 
PageRank algoritm used by Google (see \cite{Langville-PageRank} and \cite{Langville-web} for expository reviews) to 
whole genome analysis based on latent semantic indexing \cite{Ng,Homayouni}.
\cut{The performance of the algorithm is analysed by use of probabilistic estimates
in \cite{Papadimitriou} and computer implementations of the methods
have been proposed in \cite{Forman}.}

From the information contained in the dataset $\sfb{X}$, a 
weighed graph of
interactions among documents  is constructed. 
To palliate the weaknesses of principal component analysis,
reproducing kernel
methods can be used. The oldest account of these methods seems to be
\cite{Mercer}
and their formulation in the context of Hilbert spaces can be found
in \cite{Aronszajn}.
In \cite{Vert}, analysis of features of a microarray experiment
is proposed based on kernel estimates on a graph. Note however
that in that paper, the graph incorporates extrinsic information
coming from participation of genes in specific pathways as documented
in the KEGG database. On the contrary, in the method we are proposing here,
the graph can be constructed in an intrinsic way, even in the absence
of any additional information from existing databases.
In \cite{Belkin,Coifman,Nilsson}, kernel methods and Laplace eigenspace
decomposition are used to generalise principal components analysis to
include non-linear interactions among genes. Particular types of
kernels, defined in terms of commuting times for a random walk on the graph are
used in \cite{Fouss-rw,Meila,Vast}. 
All these methods, although not always explicitly stated in these articles,
are as a matter of fact very closely related since the kernels, 
the weighed graph
Laplacian and  the simple random walk on the graph can be described
in a unified formalism \cite{Campanino-Petritis,Chung,Cvetkovic,Godsil,Mohar}.
It is worth noting that analysis of Laplacian of the graph
is used in many different contexts, ranging from 
biological applications
(proteins conformation \cite{proteins}, gene arrays \cite{Ng}) 
to web search \cite{Baldi-web} or image analysis \cite{ShiMalik}.

Fuzzy clustering has been introduced in \cite{Bezdek}; 
lately 
it was shown  \cite{Wang}
equivalent to probabilistic clustering if the objective function
is expressed in terms of the R\'enyi entropy.

The idea of describing the data in terms of abstract Hilbert spaces
has been used (in the context of database search) in 
\cite{Baeza-Yates,Dominich,Gardenfors,vanRijsbergen,Widdows}. 

The semantic distillation algorithm is based on a quantum-inspired
subspace projection, strongly reminiscent of the quantum
procedure of measurement. Although fully implemented
on classical computers, it shares with general quantum algorithms
features of non-distributive quantum logic 
\cite{petritis-iq,Redei}. The semantic approach of Quantum Mechanics
can be found in \cite{Redei,dallaChiara}. It is worth underlying that
the full fledged fuzzy logic induced by quantum semantics is not equivalent
to the standard fuzzy logic introduced in \cite{Zadeh}; it represents a genuine
extension of it \cite{dallaChiara}.

\section{Perspectives}

Various data sets (not only biological) 
are presently semantically distilled and the method compared with more
traditional approaches. Preliminary results obtained so far seem to confirm
the power of the method.

Several directions are in progress:
\begin{itemize}
\item Although the method is quantum-inspired, the fuzzy logic induced
is still standard fuzzy logic.  We are currently 
working on the extension to generalised
fuzzy logic induced by full-fledged quantum semantics.
\item
The graph analysis we performed provided us with degrees of specificities of
every gene in a particular context. These data can be 
reincorporated to the graph as internal degrees
of freedom of  a multi-layered graph that can be further analysed.
\item
The connections of the 
algorithm of semantic distillation with the algorithm of purification of quantum
states \cite{Maassen} introduced in the context of quantum computing are currently explored. 
\end{itemize}

 \bibliographystyle{plain}
		 \bibliography{bibliomath,bibliobio}

			    \end{document}